\documentclass[twoside,12pt, leqno]{amsart}
\usepackage{amsmath,amscd,amsthm,amssymb,amsxtra,latexsym,epsfig,epic,graphics}
\usepackage[matrix,arrow,curve]{xy}
\usepackage{graphicx}
\usepackage{tikz}  
\usepackage{color} 
\usetikzlibrary{arrows,calc} 
\usetikzlibrary{decorations.pathmorphing}
\usetikzlibrary{decorations.markings} 
\usetikzlibrary{decorations.pathreplacing} 
\usetikzlibrary{plothandlers}
\usepackage[alphabetic,lite]{amsrefs} 

\def\antiddot{\mathinner{\mkern1mu\raise1pt\vbox{\kern7pt\hbox{.}}\mkern2mu
        \raise4pt\hbox{.}\mkern2mu\raise7pt\hbox{.}\mkern1mu}}

\newcommand{\PP}{{\mathbb P}}

\newcommand{\ZZ}{{\mathbb Z}}

\newcommand{\coker}{{\rm{coker}\,}}

\newcommand{\s}{\mathcal}

\newcommand{\sB}{{\s B}}

\newcommand{\sF}{{\s F}}

\newcommand{\sO}{{\s O}}

\newcommand{\cO}{{\s O}}

\newcommand{\punkt}{\hspace{-.3ex}\raise.15ex\hbox to1ex{\Huge.}}

\def \fix#1 {{\hfill\break \bf (( #1 ))\hfill\break}}

\DeclareMathOperator{\Sym}{Sym}

\newtheorem{theorem}{Theorem}[section]

\theoremstyle{definition}

\newtheorem{remark}[theorem]{Remark}

\newtheorem{example}[theorem]{Example}

\def\PP{{\mathbb P}}

\def\cornerT#1{{T_{\Rsh \kern -1pt #1}}}

\def\fix#1{{\bf ***Fix:} #1 {\bf ***}}

\def\bT{{\bf T}}
\def\bU{{\bf U}}

\makeatletter
\def\blfootnote{\xdef\@thefnmark{}\@footnotetext}
\makeatother

\def\lbracket{{[\kern-1.5pt[}}
\def\rbracket{{]\kern-1.5pt]}}

\makeatletter
\def\Ddots{\mathinner{\mkern1mu\raise\p@
\vbox{\kern7\p@\hbox{.}}\mkern2mu
\raise4\p@\hbox{.}\mkern2mu\raise7\p@\hbox{.}\mkern1mu}}
\makeatother

\usepackage{times}

\usepackage{color}

\usepackage[breaklinks,bookmarksopen,bookmarksnumbered,urlcolor=blue]{hyperref}
\hypersetup{colorlinks=true,backref=true,citecolor=blue}

\author[David Eisenbud]{David Eisenbud}
\address{Department of Mathematics, University of California at Berkeley and the Mathematical
Sciences Research Institute, Berkeley, CA 94720, USA}
\email{de@msri.org}

\author{Daniel Erman}
\address{Department of Mathematics, University of Wisconsin,
  Madison, Wisconsin, 53706, USA}
\email{derman@math.wisc.edu}

\author[Frank-Olaf Schreyer]{Frank-Olaf Schreyer}
\address{Fachbereich Mathematik, Universit\"at des Saarlandes, Campus E2 4, D-66123 Saar\-br\"ucken, Germany}
\email{schreyer@math.uni-sb.de}

\title[Tate Resolutions on Products of Projective Spaces]{Tate Resolutions on Products of Projective Spaces: \\ Cohomology and Direct Image Complexes}
\begin{document}

\begin{abstract}
We describe the  Macaulay2 package TateOnProducts and its capabilities, which include computing cohomology tables of sheaves
on products of projective spaces and the derived category pushForward of a sheaf under a morphism from a projective scheme to a projective space.
\end{abstract}

\maketitle

\blfootnote{
Keywords: cohomology groups, Tate resolutions, Beilinson monads \smallbreak
The first and second authors are grateful to the
National Science Foundation for partial support. 
This work is a contribution to Project I.6 of the third author within the SFB-TRR 195 "Symbolic Tools in Mathematics and their Application" of the German Research Foundation (DFG).}

\section*{Introduction}
The main functions implemented in this package are:
\begin{enumerate}
	\item  {\tt cohomologyHashTable}, which computes all of the cohomology groups of a coherent sheaf on a product
	of projective spaces, within a specified range of multidegrees;
	\item  {\tt beilinsonMonad}, which computes the Beilinson monad of a sheaf on a product of projective spaces; and
	\item {\tt directImageComplex}, which computes the derived push forward of a coherent sheaf with respect to any of the natural projection maps from a product of projective spaces.
\end{enumerate}

The algorithms employed exploit the Koszul duality between polynomial rings and exterior algebras, thus turning questions of cohomology of sheaves into questions about free resolutions of modules over an exterior algebra.  The core algorithm of the package is the computation of a finite part of the {\bf Tate resolution} of a coherent sheaf, which is a doubly infinite, multigraded, free complex of free modules over an exterior algebra whose.  The functions mentioned above are easy to compute once this is done.

Tate resolutions were introduced for sheaves on $\PP^n$ in \cite{EFS}.  The Macaulay2 package {\tt BGG} implements an algorithm for computing Tate resolutions for such sheaves, and for using them to compute sheaf cohomology and Beilinson monads~\cite{M2BGG}.  Tate resolutions for products of projective spaces are introduced in \cite{EES}.  
For a product of two or more projective spaces, there is a major new difficulty, as each term of the Tate resolution is a free module of infinite rank.  
In this package we compute a part of this infinite object that we call a {\bf corner complex}.  Corner complexes play an essential role in \cite{EES}, and they are sufficient for our applications.

This paper is organized as follows.  In Section~\ref{sec:background}, we briefly review the BGG correspondence and Tate resolution.  In Section~\ref{sec:tate resolutions}, we describe the finite part of a Tate resolution that we can actually compute.  In Section~\ref{sec:cohom tables}, we discuss the command {\tt cohomologyTable}.  
In Section~\ref{sec:beilinson monad}, we provide a brief discussion of Beilinson monads, and discuss the command {\tt beilinsonMonad}.  In Section~\ref{sec:push forward}, we discuss derived push forward algorithms. In the final Section~\ref{sec:cornercomplex}, we introduce the corner complexes, which play a mayor role in the construction of the Tate resolutions and in many of our applications.

\section*{Acknowledgments}
We thank Mike Stillman for major help with aspects of this package.

\section{The BGG correspondence and Tate resolutions}\label{sec:background}

\subsection{Review of Tate Resolutions for Projective spaces}
Fix a field $k$ and an $(n+1)$-dimensional vector space $W$.  We let $S=\Sym W$ be the symmetric algebra, with generators in degree $1$.  We also set $V=W^*$ and let $E=\Lambda V$ be the exterior algebra on $E$, with generators in degree\footnote{Our conventions in this paper are consistent with the conventions in the Macaulay2 packages {\tt BGG} and {\tt TateOnProducts}.  See Remark~\ref{rmk:conventions} for a detailed comparison of how these conventions differ from those in \cite{EFS,EES}.}  $1$.   The BGG correspondence stems from two simple observations.

First, if $A, B, C$ are finite-dimensional vector spaces over $k$, then there is a natural bijection between  homomorphisms $A\otimes_kB\to C$ and homomorphisms
$ B \to C\otimes_k A^*$.

Second, if $M = \oplus_{i\geq i_0} M_i$ is a finitely generated graded $S$-module, then the module structure induces a sequence of map $W\otimes M_i \to M_{i+1}$.  By the correspondence above, this yields a sequence of maps $M_i\to M_{i+1}\otimes V$, which induces a sequence of
linear maps of graded free $E$-modules 
$$ 
bgg(M): \cdots \to M_i\otimes E(i) \to M_{i+1}\otimes E(i+1) \to \cdots
$$
A formal computation confirms that the conditions of commutativity and associativity of the action of $S$ on $M$ exactly correspond to the condition that $bgg(M)$ is a complex; that is, consecutive maps compose to 0.

In exactly the same way, a graded $E$-module $N$ gives rise to a linear  complex of free $S$-modules
$$bgg(N): \ldots \to N_i \otimes S(i) \to N_{i+1} \otimes  S(i+1) \to \ldots $$

Tate resolutions, as defined below, are exact complexes, and our computation depends on the fact that such a complex is
determined by any differential. However, in general, only a truncation of the complex $bgg(M)$ is exact. This corresponds to a truncation of $M$ at a degree sufficiently high to have a linear resolution over $S$. The following useful criterion shows that such a truncation is necessary, and sufficient:

\begin{theorem}[Reciprocity] \cite[Theorem 3.7 and Corollary 2.4] {EFS}\label{Reciprocity}. With notation as above:
\begin{enumerate}
 \item $bgg(M)$ is an injective resolution of the $E$-module $N$ if and only if
$bgg(N)$ is a free resolution of $M$.

\item These conditions are satisfied if and only if the Castelnuovo-Mumford regularity of $M$ is 0 and $M$ has no submodule of finite length.
\end{enumerate}
\end{theorem}

In their 1978 paper, Bernstein-Gelfand-Gelfand \cite{BGG} used the BGG correspondence to identify  the derived category $D^b(\PP^n)$ of bounded complexes of coherent sheaves on $\PP^n$ with the stable category of $E$-modules. In the same volume \cite{beilinson} Beilinson described
$D^b(\PP^n)$ in terms of exceptional sequences.  (See Section~\ref{sec:beilinson monad} for more on Beilinson's work.)

Tate resolutions, as introduced in~\cite{EFS}, connect the approaches of \cite{BGG} and \cite{beilinson} and provide the foundation for the package~\cite{M2BGG}.
Let $M=\oplus M_d$ be a graded $S$-module representing the coherent sheaf $\sF$ on $\PP^n$. For any $r$, the truncation $M_{\ge r} = \oplus_{d \ge r}M_d$ represents the same sheaf.  If $r\gg 0$, then Theorem~\ref{Reciprocity} shows that the  complex
$$
bgg(M_{\ge r}) : M_r\otimes E(r) \to M_{r+1}\otimes E(r+1) \to \ldots
$$ 
is acyclic. By combining this injective resolution with a minimal free resolution of $P=\ker(M_r\otimes E(r) \to M_{r+1}\otimes E(r+1))$, we obtain an infinite exact complex
$$
\bT(\sF):   \ldots \to T^{r-2} \to T^{r-1} \to T^r \to T^{r+1} \to \ldots
$$
of free $E$-modules, which depends only on $\sF$.  This is called the {\bf Tate resolution of $\sF$}. Any finite part can be found by computing syzygies over the exterior algebra.

By \cite[Theorem 4.1] {EFS}, the Betti numbers of $\bT(\sF)$ encode the ranks of the sheaf cohomology groups of $\sF$ via the following formula:
$$
\bT^d(\sF) =T^d= \sum_{i=0}^n H^i(\PP^n,\sF(d-i)) \otimes E(d-i).
$$
and thus the syzygy computation yields the dimensions of the vector spaces $H^i(\PP^n,\sF(d))$ in any bounded range. The Beilinson monad of $\sF$ can also be obtained from the Tate resolution---see Section~\ref{sec:beilinson monad}.

\subsection{Tate Resolutions for Products of Projective Spaces}\label{subsec:products}
We now consider a product of projective spaces $\PP = \PP^{n_1}\times \ldots \times \PP^{n_{t}}$.  The Cox ring of $\PP$ is a $\ZZ^{t}$-graded polynomial ring $S=\Sym (W_1\oplus W_2 \oplus \cdots \oplus W_t)$ where $\dim W_i=n_i+1$ and where the elements of $W_1$ have degree $(1,0,\dots,0)$, the elements of $W_2$ have degree $(0,1,0,\dots,0)$ and so on.  As before, we let $V_i=W_i^*$ for all $i$, and define a $\ZZ^{t}$-graded exterior algebra $E=\Lambda (V_1\oplus V_2\oplus \cdots \oplus V_t)$ where $\deg(W_i)=\deg(V_i)$.  When comparing multidegrees in $\ZZ^t$, we always take the termwise partial order.  

By a nearly identical computation as before, the BGG correspondence $M\mapsto bgg(M)$ sends a multigraded $S$-module to a linear, mulitgraded, free, $t$-fold multil-complex of $E$-modules.  For instance, if $t=2$, then $bgg(M_{\ge(i,j)})$ is the total complex of a double complex:
\begin{equation}\label{eqn:uppercorner}
\xymatrix{ \vdots & \vdots & \\
 M_{i,j+1}\otimes E(i,j+1) \ar[r]\ar[u]& M_{i+1,j+1}\otimes E(i+1,j+1)\ar[u] \ar[r] & \cdots \\
 M_{i,j}\otimes E(i,j) \ar[r]\ar[u]& M_{i+1,j}\otimes E(i+1,j) \ar[u]  \ar[r] & \cdots \\
}
\end{equation}
In general, $bgg(M)$ is the total complex of a $t$-fold complex.

\begin{remark}\label{rmk:conventions}
We follow the conventions of the package \cite{TateOnProducts}, which differs from the conventions of \cite{EES}. In particular, in the paper:
\begin{itemize}
\item The exterior algebra $E$ is negatively graded.
\item We use $\omega_E$ instead of $E.$
\item Tate resolutions are cochain complexes instead of chain complexes.
\end{itemize}
\end{remark}
The theory of Tate resolutions for 
$
\PP$ is developed in~\cite{EES}.  If $M$ is a graded $S$-module representing $\sF$, then the Tate resolution $\bT(\sF)$ is an
exact complex of free $E$-modules with terms
\begin{equation}\label{eqn:Tate products}
\bT^d(\sF)= \sum_{a \in \ZZ^t} H^{d-|a|}(\PP,\sF(a)) \otimes E(a),
\end{equation}
where $|a|=a_1+\ldots+a_n$ denotes the total degree.
This shows that the Tate resolution encodes all the cohomology groups
$ H^j(\PP, \sF(a))$.
However, \eqref{eqn:Tate products} also shows that, when $t>1$, each term $\bT^d(\sF)$ of the Tate resolution will be an infinitely generated $E$-module. Computing the entire Tate resolution $\bT^d(\sF)$ is thus infeasible.  We can nevertheless  effectively compute the subquotient complex of $\bT^d(\sF)$ which consists of all terms in a finite range
$
low \le a \le high
$
of multidegrees with $low, \, high \in \ZZ^t $ and this is sufficient for all of our applications.

\section{Tate Resolutions}\label{sec:tate resolutions}
Each of the major applications in this package stems from the computation of a bounded part of the Tate resolution.
Let $M$ be a finitely generated $\ZZ^t$-graded $S$-module, and let $\sF$ be the 
corresponding sheaf on $\PP$. Let
$low\leq high\in \ZZ^t$ be multi-degrees defining an interval.  The function
\begin{verbatim}
tateResolution(M,low,high)
\end{verbatim}
computes a subquotient complex of $\bT(\sF)$ that contains all summands generated in degrees $low\leq a \leq  high$.  To compute this, Macaulay2 first uses  {\tt coarseMultigradedRegularity(M)} to find some $b\gg 0$ such that the total complex of $bgg(M_{\geq b})$ is acyclic.\footnote{This is closely related to the question of computing the multigraded regularity of $M$. Finding an efficient algorithm for that computation remains an interesting question.}  Then, it resolves the kernel of  the first map in $bgg(M_{\geq b})$, continuing until it covers all of the degrees in the desired range.  

For instance, if $t=2$ then we would resolve the kernel of the map:
\[
M_b \otimes E(b) \to 
\begin{matrix}
M_{b+(1,0)} \otimes E(b+(1,0))\\
\bigoplus
\\
M_{b+(0,1)} \otimes E(b+(0,1)).
\end{matrix}
\]
This yields a multigraded complex of finitely generated, free $E$-modules which is denoted $\operatorname{Tail}_b(M)$ in~\cite[\S1]{EES}.  A key fact from ~\cite[\S1]{EES} is that $\operatorname{Tail}_b(M)$ equals the subquotient complex of $\bT(\sF)$ obtained by restricting to degrees $\leq b - (1,1,\dots,1)$.\footnote{This coincidence is essential to the computability of $\bT(\sF)$ within the desired range.  It relies on properties of corner complexes, briefly described in Section \ref{sec:cornercomplex} and treated in detail in~\cite[\S3]{EES}.}
 By increasing $b$ and/or by computing enough terms of the free resolution $\operatorname{Tail}_b(M)$, we can thus compute any portion of $\bT(\sF)$.  

For example, let $\PP=\PP^1\times \PP^2$ and let $M=S^1$ be a module representing $\sF=\cO_{\PP^1\times \PP^2}$.  The input
\begin{verbatim}
T = tateResolution(M,{-3,-3},{0,0})
\end{verbatim}
computes a multi-graded complex of the form
{\small \begin{verbatim}
         1      2      6      14      29      55      97
 0  <-- E  <-- E  <-- E  <-- E   <-- E   <-- E   <-- E 
                                                                  
 -1     0      1      2      3       4       5       6  
\end{verbatim}}
\noindent which includes all the free summands in the Tate resolution that have generators in degree $(-3,-3)\leq a\leq (0,0)$. The numbers below the free modules are the homological degrees.
The Betti table of this complex with respect to total degree, is
\goodbreak
{\small \begin{verbatim}
betti T
        0 1 2  3  4  5  6
 total: 1 2 6 14 29 55 97
     0: 1 . .  .  .  .  .
     1: . 1 2  3  4  5  6
     2: . 1 3  6 10 15 21
     3: . . 1  5 15 35 70
\end{verbatim}}

\section{Computing Cohomology Tables}\label{sec:cohom tables}
To compactly encode the dimensions of the sheaf cohomology groups of a coherent sheaf $\sF$, we introduce the {\bf Euler polynomial}
$$
\sum_{i\geq 0} (\dim H^i(\PP, \sF))h^i\in \ZZ[h].
$$
For concision, we often write $H^i(\sF)$ in place of $H^i(\PP, \sF)$.
Let $M$ be a finitely generated $\ZZ^t$-graded $S$-module, and let $\sF$ be the 
corresponding sheaf on $\PP$. Let
$low\leq high\in \ZZ^t$ be multi-degrees defining an interval. We can compute
all the cohomology vector spaces of all the twists $\sF(a)$ for $low\leq a\leq high$ with the function 
{\small \begin{verbatim}
eulerPolynomialTable(M,low,high)
\end{verbatim}}
The output of this function is a hash table consisting of the pairs
$$ 
a \Rightarrow \sum_{i\geq 0} (\dim H^i(\PP, \sF(a))h^i \ \  \text{ where }\ \  a\in \ZZ^t.
$$
As in the case of a single projective space, these are computed as Betti numbers of the appropriate finite  part of the Tate resolution.

For example,  if $\PP=\PP^1\times \PP^2$ and   $\sF = \sO_{\PP^1\times \PP^2}$ is represented by $S^1$, then:
{\small \begin{verbatim}
 (S,E) = productOfProjectiveSpaces{1,2};
 low  = {-3,-3};high = {3,3};
 eT = eulerPolynomialTable(S^1, low,high);
 eT#{2,-3}
\end{verbatim}}
\noindent gives output
\begin{verb}
3h2
\end{verb}.
This means that $H^2(\sF(2,-3))$ has rank $3$ and that $H^i(\sF(2,-3)) = 0$ for $i\ne 2$.

In the case $t=2$, we can display the hash table as a matrix using the function
\begin{verb}
 cohomologyMatrix(M,low,high)
\end{verb},
where the upper right hand corner of the table corresponds to the multi-index $high$: Continuing with the previous example,
we have:
{\small \begin{verbatim}
cohomologyMatrix(S^1, low,high)
\end{verbatim}}
\noindent gives output
{\small 
\begin{verbatim}
     | 20h 10h 0 10 20  30  40  |
     | 12h 6h  0 6  12  18  24  |
     | 6h  3h  0 3  6   9   12  |
     | 2h  h   0 1  2   3   4   |
     | 0   0   0 0  0   0   0   |
     | 0   0   0 0  0   0   0   |
     | 2h3 h3  0 h2 2h2 3h2 4h2 |
\end{verbatim}}
\noindent Since $high = \{3,3\}$, the index of the top row is 3 as is the index of the right-hand column. 
The $(i,j)$-th entry is the Euler polynomial of $\sF(i,j)$. 

\section{The Beilinson Monad}\label{sec:beilinson monad}
Let $\sF$ be a coherent sheaf on $\PP$.  The Beilinson monad of $\sF$ is a complex of sheaves $\sB$ whose terms are direct sums of sheaves from an exceptional collection on $\PP$ built from exterior powers of the cotangent bundles on the factors, and whose homology groups satisfy
\[
H^{i}\sB = \begin{cases} \sF & \text{ if } i=0\\
0& \text{ if } i \ne 0 
\end{cases}.
\]
The advantage of this representation is that the terms of the $\sB$ depend only on the cohomology groups of
$\sF$. See~\cite{beilinson,AO,EFS,Huy} for discussion of Beilinson monads and ~\cite[\S2]{EES} for a discussion in the context of $\PP$.  As explained in~\cite[\S2]{EES}, the Beilnson monad for $\sF$ is determined by the subquotient complex of $\bT(\sF)$ consisting of summands generated in the finite range of degrees $0\leq a \leq n$.  This range of degrees is known as the {\em Beilinson window,} and a Beilinson monad for a $\sF$ can be  computed easily from any finite part of the Tate resolution that contains the Beilinson window.  

For the computation, we start with a module $M$ representing a coherent sheaf $\sF$ on $\PP$. Executing
{\small \begin{verbatim}
B = beilinsonMonad M
\end{verbatim}}
\noindent produces a complex $B$ of graded $S$-modules such that the corresponding complex $\widetilde{B}$ of sheaves on $\PP$ is the Beilinson monad for $\sF$.
    
For example, if we take $M$ to be a twist of the syzygy module of the maximal ideal in the bihomogeneous
coordinate ring $S$ of $\PP=\PP^1\times \PP^2$, we can execute:
{\small \begin{verbatim}
M =  S^{{1,1}} ** ker vars S;
\end{verbatim}}
The command \verb beilinsonWindow \ computes the subquotient complex of the Tate resolution in the degrees $0\leq a \leq n$.  In our example, this complex only has two terms:
{\small \begin{verbatim}
T=tateResolution(M,low,high);
W = beilinsonWindow T
        6      1
       E  <-- E
betti W
             0 1
      total: 6 1
          0: 6 .
          1: . 1
 \end{verbatim}}
The Beilinson monad itself is obtained by applying a certain $\bU$ functor (see \cite[\S2]{EES}) to any piece of $\bT(\sF)$ which contains the Beilinson window.  This whole process is computed directly from the module $M$ via
\begin{verbatim}
B = beilinsonMonad M
\end{verbatim}
which, in this case, gives a two-term complex:
{\small
\begin{verbatim}     
 6
S  <-- cokernel {1, 1} | x_(1,2)  |
                {1, 1} | -x_(1,1) |
                {1, 1} | x_(1,0)  |
\end{verbatim}}
Macaulay2 computes sufficiently many terms of the Tate resolution
of $\sF$ (as in \S\ref{sec:tate resolutions}), extracts the subquotient complex corresponding to the Beilinson window, and then applies the functor $\bU$.

The 0th-cohomology of the Beilinson monad is a module representing the sheaf $\sF$,
while the other cohomology groups are trivial as sheaves---that is, supported on the irrelevant ideal. 
Note also that the 0th-cohomology may not be equal to the module $M$, but the two modules
will agree after  truncating at a sufficiently positive multi-degree.
For instance:
{\small
\begin{verbatim}
isIsomorphic(HH^0 B ,M)
   false
isIsomorphic (truncate({0,0},HH^0 B),truncate({0,0},M))
   true
\end{verbatim}}
\noindent Similarly, the other cohomology of $B$ will have irrelevant support, and therefore will be 0 in sufficiently positive\footnote{For the precise meaning of a {\bf sufficiently positive degree} see~\cite[Definition~1.8]{EES}.} multi-degrees.

\section{The direct image complex}\label{sec:push forward}
This package can also compute the pushforward complex of a sheaf along a morphism
to projective space. The main ingredient is the pushforward to a factor in a product of projective
spaces, which we explain first. 

Let $I \subset \{0,\ldots,t-1 \}$ be a proper subset and let
$$ \pi \colon\PP \to  \PP^I:=\prod_{i \in  I} \PP^{n_{i+1}}$$
be the projection to the corresponding factors. (Note that Macaulay2 prefers integers sequences starting with $0$ in many commands, which is the reason for the re-indexing.) Given a coherent sheaf $\sF$ on $\PP$, represented by an $S$-module $M$,
the function
{\small \begin{verbatim}
directImageComplex(M,I)
\end{verbatim}}
\noindent computes a complex of $S_I = \Sym (\oplus_{i \in I} W_{i+1})$-modules, whose sheafification
is a Beilinson monad for the direct image complex $R\pi_* \sF$ in $D^b(\PP^I)$.   See~\cite[Corollary 0.3]{EES} for the related theorem.

The key idea is that an appropriate multi-graded strand of the Tate resolution of $\sF$ is a Tate resolution for $R\pi_*(\sF)$ on the image.  (While we have focused on Tate resolutions of sheaves, one can also define the Tate resolution of an element in the derived category, such as $R\pi_*(\sF)$.)

For example, let $\sF=\cO_{\PP^1\times \PP^2}$, represented by the module $S^1$ and let $\pi:\PP^1\times \PP^2\to \PP^1$.  As in \S3, we have:
{\small 
\begin{verbatim}
cohomologyMatrix(S^1, {-3,-3},{3,3})
     | 20h 10h 0 10 20  30  40  |
     | 12h 6h  0 6  12  18  24  |
     | 6h  3h  0 3  6   9   12  |
     | 2h  h   0 1  2   3   4   |
     | 0   0   0 0  0   0   0   |
     | 0   0   0 0  0   0   0   |
     | 2h3 h3  0 h2 2h2 3h2 4h2 |
\end{verbatim}
} 
The far left vertical column of this matrix corresponds to the subquotient of $\bT(\sF)$ consisting of summands generated in degree $(-3,i)$ for some $-3 \leq i\leq 3$.  That subcomplex, which only involves exterior variables corresponding to the $\PP^2$-factor, is a finite part of the Tate resolution of $R\pi_*(\sF(-3,0))$ on $\PP^2$.  It can be computed by
\begin{verbatim}
T = tateResolution(S^1,{-3,-3},{3,3});
s = strand(T,{-3,0},{0});
\end{verbatim}
Note that by the K\"unneth formula, $R^0\pi_* \sF(-3,0))=0$ and
\[
R^1\pi_*\sF(-3,0)\cong H^1(\PP^1,\cO_{\PP^1}(-3))\otimes \cO_{\PP^2}[-1] \cong \cO_{\PP^2}^2[-1]
\]
where $[-1]$ indicates that the sheaf lies in cohomological degree $1$.   We can view the Betti table of this strand:
\begin{verbatim}
betti s
             -1  0 1 2 3 4  5  6
o76 = total: 20 12 6 2 2 6 12 20
          1: 20 12 6 2 . .  .  .
          2:  .  . . . . .  .  .
          3:  .  . . . 2 6 12 20
\end{verbatim}
Note that, since it is a strand, the terms retain the same homological that they had within the Tate resolution.  Combining this with the methods of \S\ref{sec:beilinson monad}, we can compute a Beilinson monad for $R\pi_*(\sF(-3,0))$ directly from this strand of its Tate resolution.

\medskip

We now show how this enables one to solve the general problem of computing $R\pi_* \sF$ of a sheaf $\sF$ on a projective scheme $X$ under a morphism 
$
\pi \colon X \to \PP^m.
$
To specify the morphism $\pi$, we suppose that $X$ is given by its homogeneous ideal $J \subset R = K[x_0,\ldots,x_n]$. Locally near any point of $X$ the morphism can be given by $m+1$ homogeneous
forms of the same degree, but a different representation is necessary at points where all these forms are 0. To represent $\pi$
at every point of $X$ at once, we thus allow $k$ such representations, which we package as a
 $k \times (m+1)$ matrix $\phi=(\phi_{ij})$ of homogeneous forms in $R$, such that any two forms
 in the same row have the same degree. To ensure that the different rows define the same
 morphism where both are nonzero, we insist that the $2\times 2$ minors of $\phi$ are contained in $J$;
 to ensure that the morphism has no base locus, 
 we insist that the entries of $\phi$ have no common zero on $X=V(J)$.

Represent $\sF$ by an $R$-module $N$
(whose annihilator contains $J$) and let $T$ denote the homogeneous coordinate ring of $\PP^m$.  The function
{\small \begin{verbatim}
directImageComplex(J,N,phi)
\end{verbatim} }
\noindent returns a complex of $T$-modules whose sheafification
is a Beilinson monad for $R\pi_* \sF$. 

\begin{example}
Consider the cubic scroll 
$\PP(\sO_{\PP^1}(2) \oplus \sO_{\PP^1}(1))\cong X\subset \PP^4
$ 
defined by the ideal of minors of the $2\times 3$ matrix
$$
m=\begin{pmatrix}
x_0& x_1 & x_3 \cr
x_1 & x_2 & x_4 \cr
\end{pmatrix}   .
$$
Let $\pi :X \to \PP^1$ be defined at each point by the ratio of (at least) one of 
the columns of $m$, so that the map $\pi$ is represented, in the sense above,
by the matrix $\phi = {\rm transpose}\ m$. 

Let $N=\Sym^2(\coker m)\otimes R(1)$.  The sheaf $\sF$ associated to $N$ is the line bundle $\sF = \sO_X(1)\otimes \pi^*\sO_{\PP^1}(2)$ on $X$.  We compute:
{\small\begin{verbatim}
kk = ZZ/101;
R = kk[x_0..x_4];
m = matrix{{x_0,x_1,x_3},{x_1,x_2,x_4}}
J = minors(2,m);
N = symmetricPower(2,coker m)**R^{1};
phi = transpose m;
RpiN = directImageComplex(J,N,phi);
T = ring RpiN;
RpiN
        9      7
       T  <-- T               
       0      1
\end{verbatim}    }   
\noindent We have thus $R\pi_*\sF$ is represented by the complex $\cO_{\PP^1}^9\gets \cO_{\PP^1}(-1)^7$.  Using this, one can check that $R^1\pi_*\sF=0$ and $R^0\pi_*\sF\cong  \sO_{\PP^1}(3) \oplus \sO_{\PP^1}(4)$.
\end{example}

\section{Corner complexes.}\label{sec:cornercomplex}
Let $M$ be a finitely generated, multi-graded module representing a coherent sheaf $\sF$ on $\PP$ and let $c$ be a multi-degree such that $ bgg(M_{\geq c})$ is acyclic.  
The corner complex is an exact complex of the form
\[
\operatorname{Tail}_c(M) \to bgg(M_{\geq c}).
\]
The exactness of the corner complexes, and the fact that $\operatorname{Tail}_c(M)$ equals a subquotient complex of $\bT(\sF)$ in appropriate degrees, are the essential features which allow one to compute Tate resolutions in the first place.  

Since corner complexes are so fundamental to this theory, we also introduced a function {\tt cornerComplex} for computing any portion of a corner complex at $c$, at least assuming that $bgg(M_{\geq c})$ is acyclic.  In \S\ref{sec:tate resolutions}, we described how to compute $\operatorname{Tail}_c(M)$ in this case.  By also computing enough terms of the injective resolution $bgg(M_{\geq c})$, which by \ref{subsec:products} amounts to linear algebra, we can thus compute any portion of a corner complex.

\begin{example}
We again consider 
$\PP^1\times \PP^2$, with homogeneous coordinate ring $S$.
We compute the corner complex of $S^1$ at $c=\{0,0\}$:
{\small 
\begin{verbatim}
r = cornerComplex(S^1,{0,0},{-4,-4},{3,3});
betti r
             -6 -5 -4 -3 -2 -1 0 1 2  3  4  5   6   7
      total: 40 54 50 35 15  5 1 1 5 15 35 70 126 210
          0: 40 54 50 35 15  5 1 . .  .  .  .   .   .
          1:  .  .  .  .  .  . . . .  .  .  .   .   .
          2:  .  .  .  .  .  . . . .  .  .  .   .   .
          3:  .  .  .  .  .  . . . .  .  .  .   .   .
          4:  .  .  .  .  .  . . 1 5 15 35 70 126 210
\end{verbatim}
}
This computation thus recovers terms corresponding to sheaf cohomology groups $H^0$ and $H^3$.  In particular, the corner map goes from $E^1(-2,-3)$ to $E^1$.  If we were to increase $c$, then we would recover some $H^1$ nad $H^2$ groups as well.
\end{example}
For further information on corner complexes, see \cite[\S3]{EES} and the documentation for the function  {\tt composedFunctions}.

\end{document}